\documentclass[11pt]{article}
\usepackage{times}
\usepackage{graphicx}
\usepackage{mathtools}
\usepackage{color}

\title{M\"obius strips before M\"obius: \\ Topological hints in ancient representations}

\author {Julyan H. E. Cartwright$^{*,\dag}$ and Diego L. Gonz\'alez$^{\ddag,\S}$ \\
\\ 
$^*$ Instituto Andaluz de Ciencias de la Tierra, \\ CSIC--Universidad de Granada, \\ E-18100 Armilla, Granada, Spain  \\
$^\dag$ Instituto Carlos I de F\'{\i}sica Te\'orica y Computacional, \\ Universidad de Granada, E-18071 Granada, Spain \\
$^\ddag$ Istituto per la Microelettronica e i Microsistemi, \\ Area della Ricerca CNR di Bologna,  I-40129 Bologna, Italy \\
$^\S$Dipartimento di Scienze Statistiche ``Paolo Fortunati'', \\ Universit\'a di Bologna, I-40126 Bologna,  Italy 
}

\date{Version 1.9 of \today}

\begin{document}

\maketitle

August M\"obius discovered his eponymous strip  --- also found almost contemporaneously by Johann Listing --- in 1858, so a pre-1858 M\"obius band would be an interesting object.  It turns out there were lots of them.

\section*{M\"obius mosaics}

\begin{figure}[tbp]
\includegraphics*[width=\textwidth]{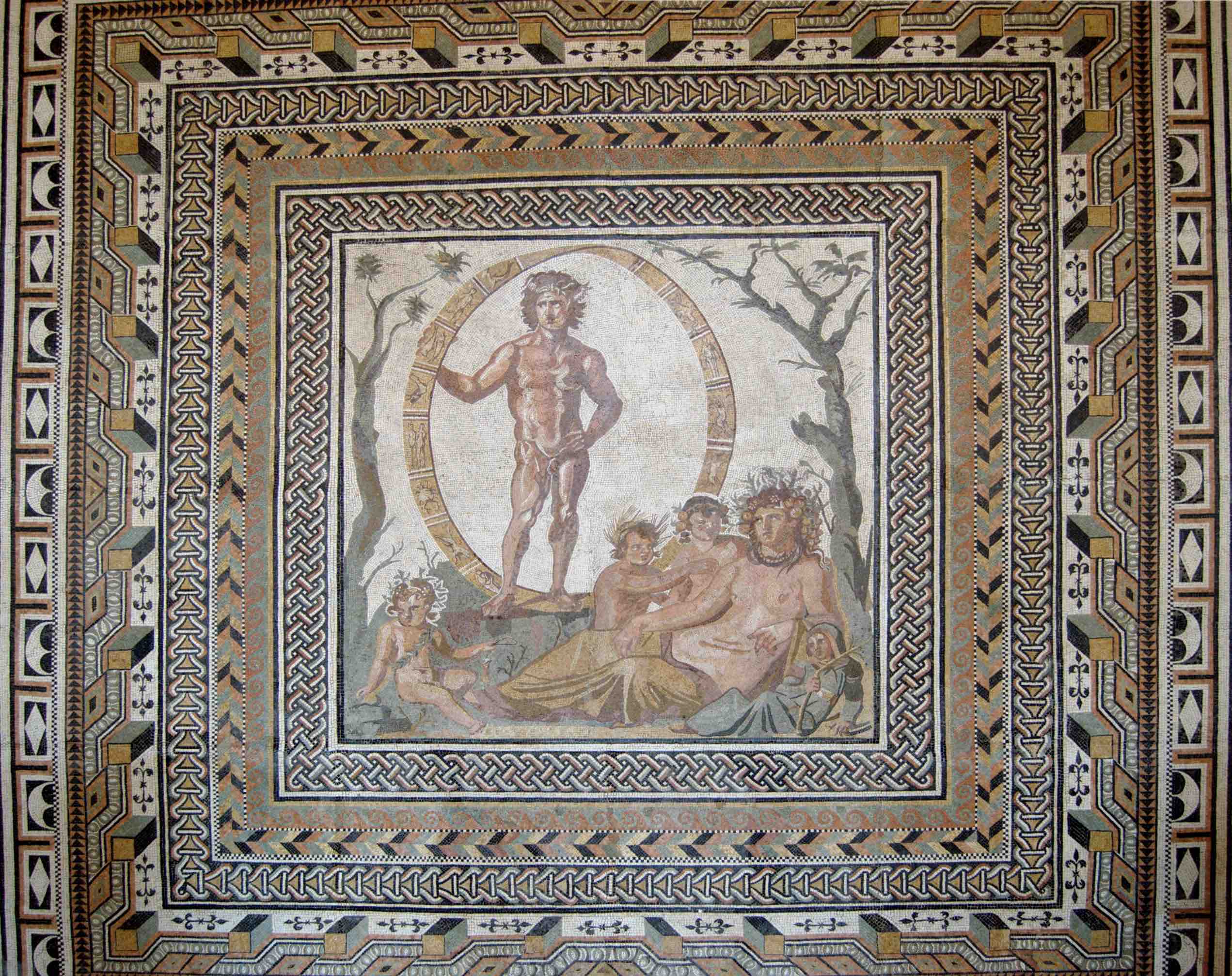}
\caption{\label{aion1}
Mosaic showing Aion and the Zodiac from a Roman villa in Sentinum, today Sassoferrato in Marche, ca.\ 200--250 CE. Aion, god of eternal time, stands in the celestial belt of the ecliptic decorated with zodiacal signs, between a green and a bare tree (summer and winter). Before him is the mother-earth Tellus, the Roman Gaia, with four children, the four seasons personified. The original is displayed at the Glyptothek Museum in M\"unich.  Image from Bibi Saint-Pol / Wikimedia Commons.
}
\end{figure}

\begin{figure}[tbp]
\includegraphics*[width=\textwidth]{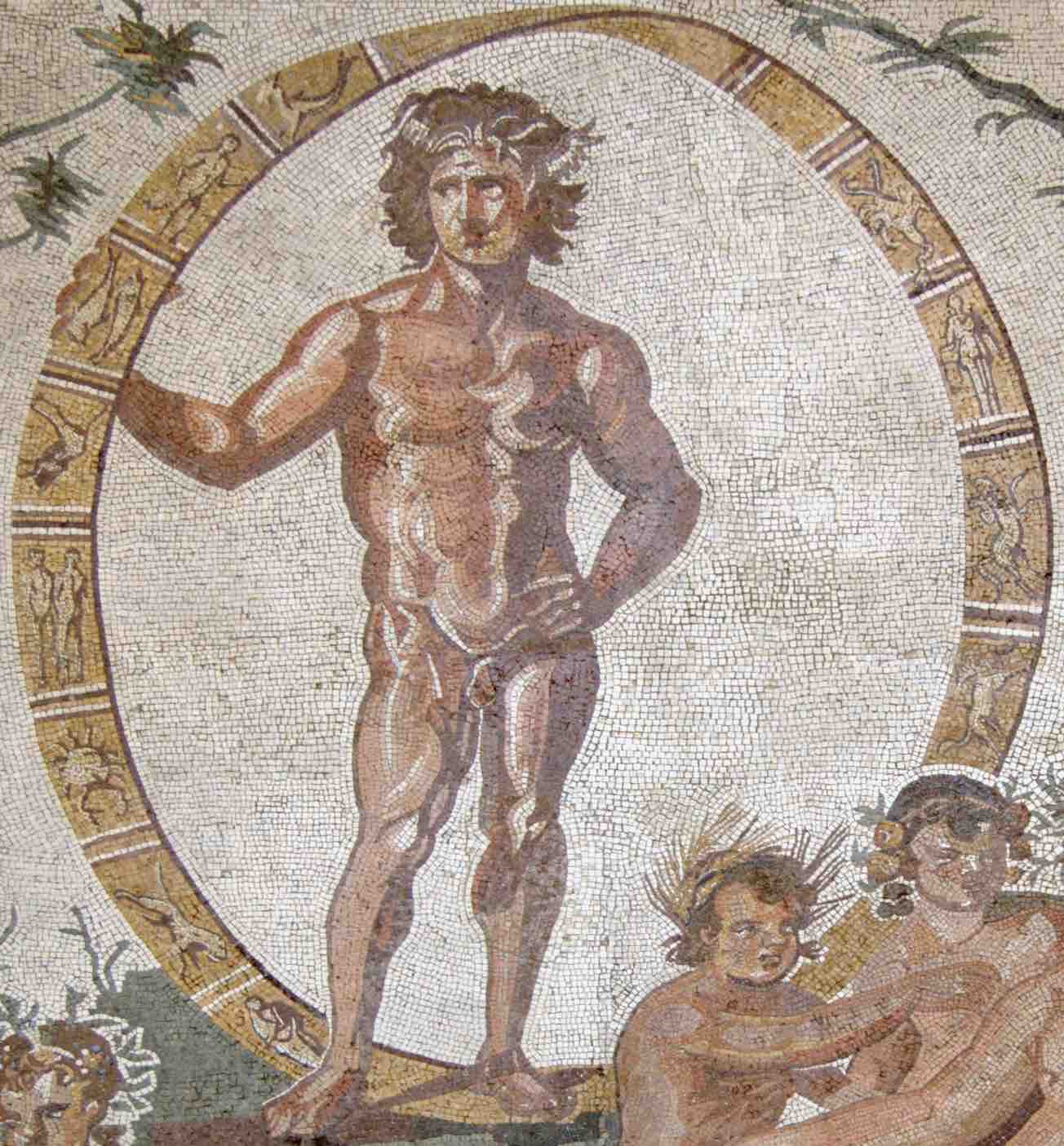}
\caption{\label{aion2}
Central part of the great Sentinum  mosaic showing the M\"obius band.  Image from Bibi Saint-Pol / Wikimedia Commons.
}
\end{figure}

At Sentinum around 295 BC there took place perhaps the greatest battle of antiquity in the Italian peninsula. The engagement marked the true beginning of the Roman empire; at it the Samnites, Umbrians, Gauls, and Etruscans were definitively defeated and submitted to Roman authority. The civic museum of Sassoferrato --- which is the present-day name of the nearby town that originated in the ferrous minerals contained in the surrounding hills --- hosts in its archaeological section the fruit of more than two centuries of  excavations, including most of the material remains of the city of Sentinum, which existed until around the sacking of Rome in 410 CE. 
Among these antiquities, there is a reproduction of a polychrome mosaic representing Aion and the zodiac; Fig.~\ref{aion1}. 

The Sentinum mosaic, dated to 200--250 CE,  depicts a M\"obius strip; Fig.~\ref{aion2}. The zodiacal symbols are shown in the band that surrounds Aion. At the bottom, Aion is standing with one foot on the internal part of the band. If this is followed clockwise, we arrive at the external part of the band over Aion's head. If, instead, we follow the band anticlockwise, we arrive at the internal part of the band at the same position over his head, which demonstrates the M\"obius topology of the strip. It is worthy of note that the zodiacal symbols of the constellations are not all in their natural places in the band of the ecliptic; was this mixing up an error on the part of the artist, or a deliberate choice with some hidden, perhaps mystical meaning \cite{Venturini2008}?

\begin{figure}[tbp]
\includegraphics*[width=\textwidth]{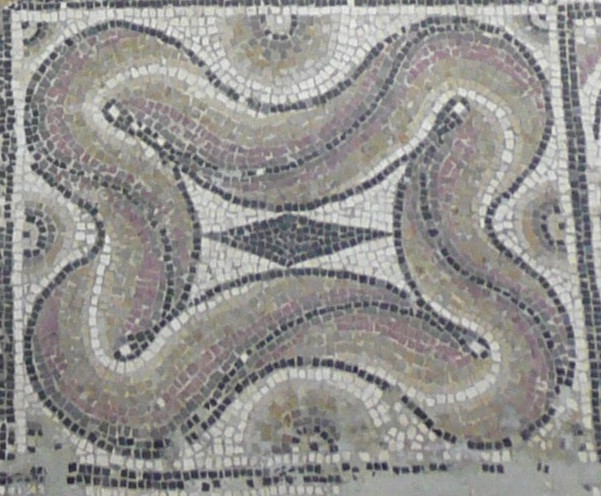}
\caption{\label{larison}
Roman mosaic representation of a M\"obius band with five half-twists discussed by Larison  \cite{Larison1973}; detail from a mosaic depicting the legend of Orpheus preserved at the Mus\^ee Lapidaire d'Art Pa\"{\i}en, Arles. 
Image from Finoskov / CC-BY-SA-3.0 / Wikimedia Commons.
}
\end{figure}

\begin{figure}[tbp]
\includegraphics*[height=0.33\textwidth]{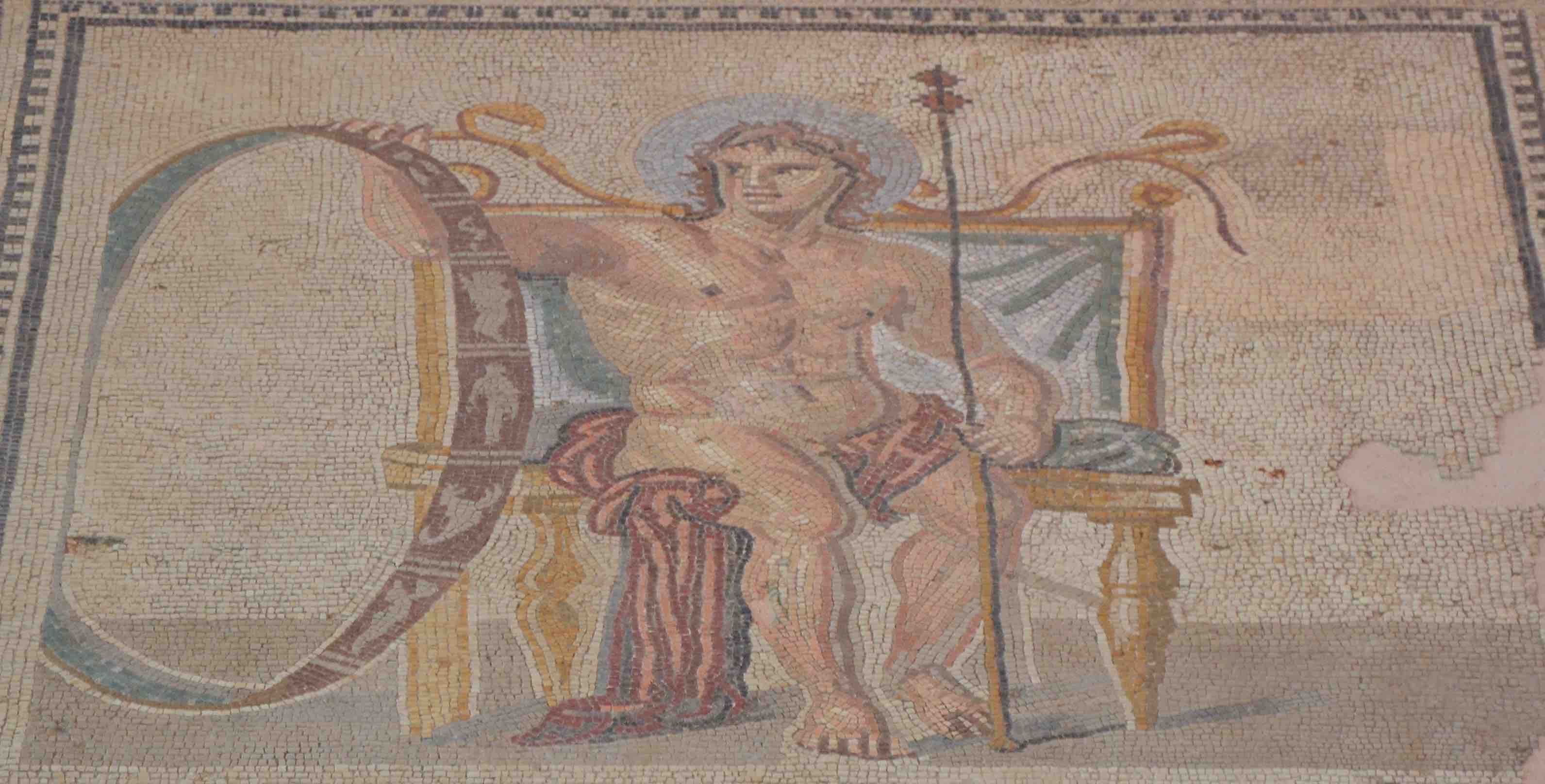}
\includegraphics*[height=0.33\textwidth]{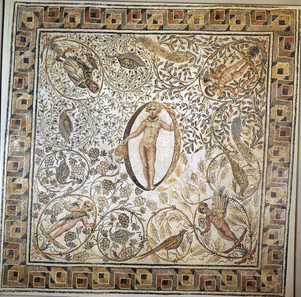}
\includegraphics*[height=0.455\textwidth]{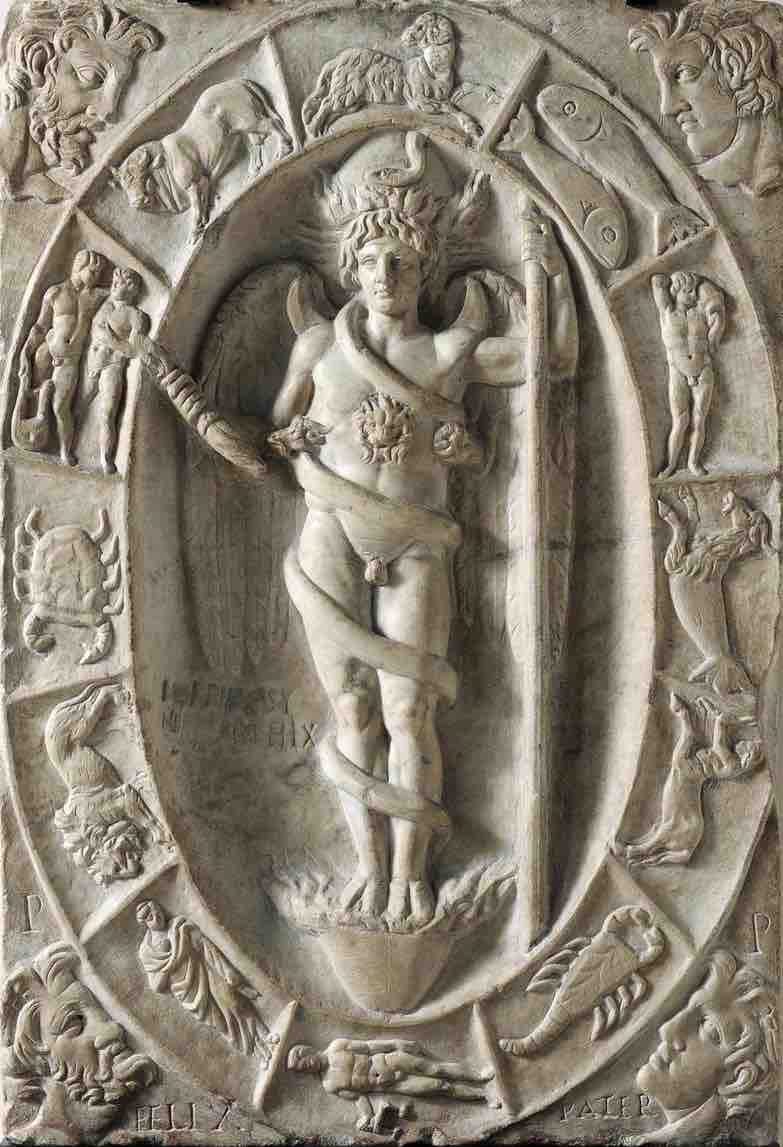}
\includegraphics*[height=0.455\textwidth]{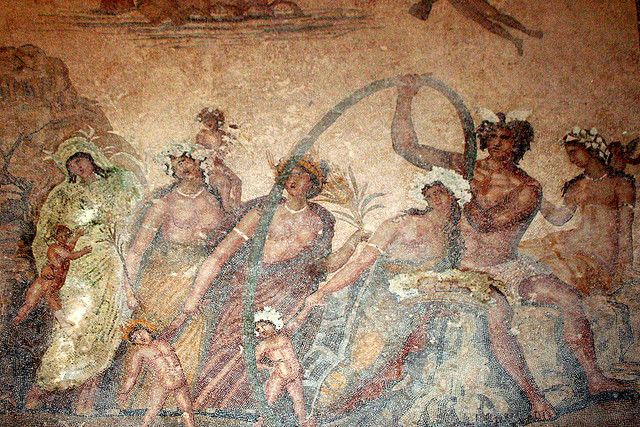}
\includegraphics*[height=0.485\textwidth]{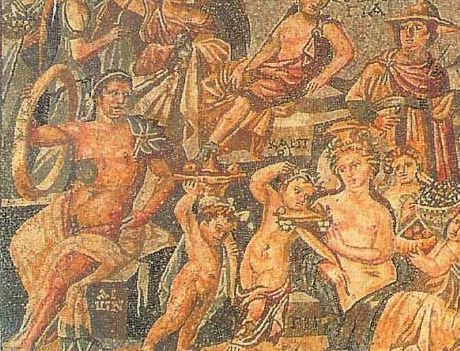}
\includegraphics*[height=0.485\textwidth]{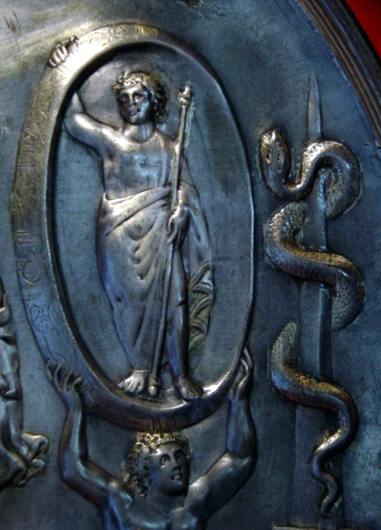}
\caption{\label{aion3}
Non-M\"obius ancient Roman representations of Aion and the zodiac:
(a) Mosaic of 3rd century CE from the triclinium of a  house in Arles, Mus\'ee de l'Arles antique, Arles. 
(b) Mosaic from 4th century CE from Ha\"{\i}dra, modern-day Tunisia, now in the United Nations building in New York \cite{Parrish1995}.
(c) Relief of 2nd century CE now in the Galleria Estense in Modena.
(d) Mosaic from Villa Selene, near Leptis Magna,  modern-day Libya.
(e) Mosaic from  Shahba--Philippopolis \cite{Quet1999}, now in Damascus Museum, Syria.
(f) Detail from the Parabiago  plate, possibly 2nd century CE, Museo Archeologico, Milan.
 Images from  Giovanni Dall'Orto, Carole Raddato  / CC-BY-SA-3.0 / Wikimedia Commons.
}
\end{figure}

The original mosaic is displayed at the Glyptothek Museum in M\"unich, where it was taken in 1828. The M\"obius strip was discovered three decades later in (what is today) Germany by both M\"obius and Listing in 1858 \cite{Fauvel1993,Pickover2006}. A fascinating, but probably unanswerable question is then whether M\"obius, and/or Listing might have seen this mosaic in the museum in M\"unich and might, perhaps unconsciously, have been influenced by it in their discoveries of the band.

The appearance of the M\"obius strip in Roman mosaics was first proposed in a very interesting article  by Lorraine Larison \cite{Larison1973}. Our example is less dubitable than those she discussed, principally one from Arles  showing a band with five half-twists (Fig.~\ref{larison}) provoked some doubts as to whether it was in fact representing a M\"obius strip \cite{Dickerson1974}. In the case of the Sentinum mosaic there is no such room for doubt.  
There are many interpretations of the god Aion, most of them related to the ideas of cyclical or unbounded time; Aion has been identified with Kronos and Uranus, among others. 
Other ancient representations of Aion holding the zodiac are known, both --- interestingly, like Larison's mosaic --- from Arles and from other parts of the Roman empire; Fig.~\ref{aion3}. However, the zodiacal hoops in these other representations  are not M\"obius strips. Of course, a circle should represent also the cyclicity of time, but a M\"obius strip, implying a doubling of the time to return to the original starting point, gives a more impressive idea. 
Regarding this point, we can remark that mathematically the invariant manifold in which a period-doubling bifurcation is born is indeed a M\"obius strip \cite{Ruelle1989}. A wonderful physical example is represented by the movement of charged particles trapped in the Earth's magnetic field that, just before their behaviour becomes chaotic, lie precisely on a M\"obius strip \cite{Ynnerman2000}.

\section*{More twisted bands from mosaics, the ouroboros, and infinity}

\begin{figure}[tbp]
\centering\includegraphics*[height=0.4\textwidth]{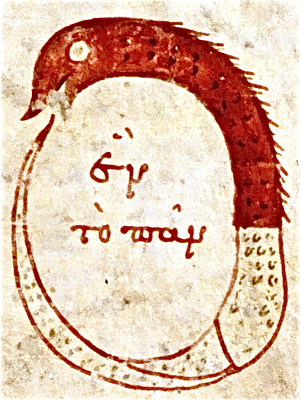}
\centering\includegraphics*[height=0.4\textwidth]{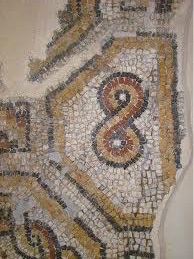}
\centering\includegraphics*[height=0.4\textwidth]{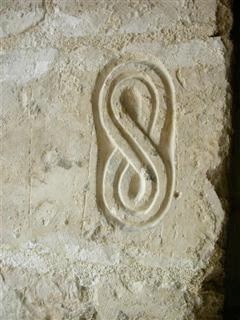}
\caption{\label{ouroboros1}
(a) An ouroboros from the  \emph{Chrysopoeia of Cleopatra}  \cite{Chrysopoeia} proposed to be a M\"obius band.
(b) Roman mosaic ``infinity'' symbol.
(c) The double-band ``infinity'' symbol at the abbey of San Vittore alle Chiuse in the Marches, Italy.
(Note the opposite chiralities of (b) and (c).)
}
\end{figure}

It has been hypothesized that the ouroboros, the ancient symbol of the snake devouring its own tail, thus symbolizing the cyclic character of time, may be an ancient representation of a M\"obius band \cite{D'amore2015}. The first known appearance of the ouroboros motif is in the \emph{Enigmatic Book of the Netherworld} \cite{netherworld}, an ancient Egyptian funerary text of the 14th century BC. An illustration from the \emph{Chrysopoeia of Cleopatra} \cite{Chrysopoeia}, almost contemporary with the Sentinum mosaic, depicts an ouroboros that has been proposed to be a M\"obius band (Fig.~\ref{ouroboros1}(a)). The \emph{Chrysopoeia of Cleopatra} is a sheet of papyrus supposed to have been written by Cleopatra the Alchemist in the 3rd century CE. However, it is not clear from a visual point of view that the representation symbolizes precisely a M\"obius band. Some assumptions about the way in which the band can be represented in a two-dimensional surface need to be made, such as, for example, the existence of lines indicating the crossing of  borders. However, a striking connection with the Sentinum mosaic is that the ouroboros is also a symbolic attribute of the god Aion.

A connection between the M\"obius band and the figure of eight lying on its side that is now known as the symbol of infinity has often been proposed.  A physical M\"obius band often adopts such a figure-of-eight form due to the internal tensions related to the half-twist topology (we may recall Escher's famous \emph{M\"obius strip II} in which ants crawl on such a figure of eight). What is beyond doubt is that John Wallis introduced the current mathematical meaning of the symbol in his 1655 work \emph{De sectionibus conicis} \cite{conicis}. It is possible that Wallis may have had in mind the symbol as a form of the Roman numeral for one thousand, which is today usually written M, but in the past was sometimes written C\reflectbox{C} or CI\reflectbox{C}, and thereby $\infty$ \cite{Menninger}.
 The same symbol --- albeit probably simply having the same pattern, rather than implying the same meaning --- is also found in many Roman mosaics; Fig.~\ref{ouroboros1}(b). 

An interesting version of the infinity sign is also found not very far from Sentinum. An engraved stone in the Romanesque abbey of San Vittore alle Chiuse in the Marches shows a three-dimensional representation of the infinity sign using two interlinked bands: Fig.~\ref{ouroboros1}(c). Each band forms separately an asymmetric figure of eight, with a major and minor loop, but their superposition restores the overall symmetry. 

The infinity signs of Figs~\ref{ouroboros1}(b) and \ref{ouroboros1}(c) are not M\"obius strips, as may readily be checked by tracing an edge; they are bands with two half-twists, or one full twist. When they are cut down the middle, `bisected' in the jargon, the outcome will be two interlinked strips each with two half-twists. A band with two half-twists cannot be obtained by bisecting a M\"obius strip, either: when a half-twist M\"obius strip is bisected, the result is a strip with not two but four half-twists, or two full twists. And with M\"obius strips of more half-twists, one obtains strips of yet more twists; the rule is that a strip with $N$ half-twists, when bisected, becomes a strip with $N + 1$ full twists. 

The relation between the infinity symbol of the type shown in our Figures~\ref{ouroboros1}(b) and (c)  and a M\"obius strip is thus indirect: both arise from twisting a belt before joining its ends, but with different numbers of half-twists. Observations show that the ancient Romans liked to play with the ideas of bands in the patterns of their mosaics, and it seems to us very likely that Roman mosaic artists may have played with physical bands and have made examples with varying numbers of half-twists to test the pattern when designing mosaics.

\section*{M\"obius drive belts}

The M\"obius strip has certainly found at least one practical use that it is easy to conceive may have predated its discovery by M\"obius. The M\"obius band has sometimes been used --- and is sometimes used today ---  as a belt to connect machinery. The earliest written reference to this industrial application of the M\"obius strip that we have been able to trace is in the pages of \emph{Scientific American} in 1871 \cite{Anon1871}; thus just post-M\"obius:
\begin{quote}
If it would not be too presumptuous in a young mechanic, I will give to J. F. K. and others if they wish it,  a rule for putting on a quarter twist belt to make it stretch alike on both edges, and do their work well, no matter what the width of the belt. The belt is to be put on in the usual way, and the ends brought together ready for lacing. Then turn one piece the opposite side (or inside) out, and lace. The belt will run, it will be found, first one side out, and then the other and will draw alike on both sides.
\end{quote}
By `quarter twist belt' the writer meant a belt to drive a piece of machinery whose pulley is at ninety degrees to the drive shaft; refer to Fig.~\ref{inman}. 

The business of whether and when to twist a drive belt to make a M\"obius strip provoked a lively correspondence in the pages of the \emph{American Machinist} in 1903 \cite{Inman1903} (see also \cite{Anon1903,Arnold1903,Fulton1903}):
\begin{quote}
Editor American Machinist. In your issue of September 10, bottom of page 1292, you quote the Wood Worker as telling of a belt â``sewed with a half twist in it, thereby bringing both sides alternately in contact with the pulleys, with the result that the belt drove better than before.'' You add, ``This will not help us much in the solution of the perpetual question as to which side of the belt to run next to the pulley.'' 

The inference from this comment would be that the belt in question was twisted for the purpose of ascertaining which side of would give the best service under identical conditions; perhaps this true, but brings to my mind several belts have seen running with a half twist in them that was put there for an entirely different purpose.

\begin{figure}[tbp]
\centering\includegraphics*[width=0.3\textwidth]{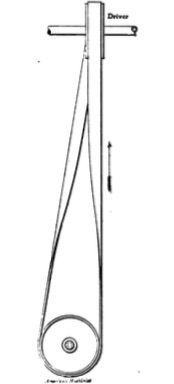}
\caption{\label{inman}
M\"obius belt drive: ``Half twisting a quarter twist belt''. From Inman \cite{Inman1903}.
}
\end{figure}

I was at one time the proprietor of a planing mill where it was necessary to transmit motion from one line shaft to another where the two were at right angles to each other;  the belt was first used in the ordinary way and ran well enough for a while, then it commenced breaking the laces and giving trouble generally. One day the belt was taken down and stretched upon the floor and it was found to be bowed --- like an iron barrel hoop when it is cut and straightened out --- a result which one can readily see the reason for upon careful observation, as the belt has more tension upon one edge than the other. The belt was then sewn together with a half  twist, as shown in the accompanying drawing [Fig.~\ref{inman}], with the result that the trouble ceased. 

Since that time I have been about the country a good deal and have found several belts running under the same conditions and giving the same trouble and have had the pleasure of suggesting the half-twist remedy, the application of which has always brought the same satisfactory results.
\end{quote}

From this correspondence we may gather that the concept of using a M\"obius belt to drive machinery was quite extended by the beginning of the 20th century. However, we have not been able to trace it back in writing to before M\"obius. The concept of introducing a twist before sewing a belt was not present in Willis' 1841 text book of mechanical engineering, \emph{Principles of Mechanism} \cite{Willis1841}, which had the following to say about belts:
\begin{quote}
Flat leather belts appear to unite cheapness with utility in the highest degree, and are at any rate by far the most universally employed of all the kinds. \ldots [they] are best joined by simply overlapping the ends and stitching them together with strips of leather passed through a range of holes prepared for the purpose \ldots Every leather belt has a smooth face and a rough face. Let the rough face be placed in contact with both pullies.
\end{quote}
We see that the utility of the natural inhomogeneity of the two sides of a leather belt is what most preoccupied Willis\footnote{Willis does mention a figure-of-eight configuration for belting, but his preferred solution has zero half-twists; see his figure 95.}; and ``which side of the belt to run next to the pulley'' was said to be the ``perpetual question'' in the \emph{American Machinist}. Clearly this inhomogeneity would be negated by a M\"obius strip version, and perhaps it is only with later materials like rubber belts that are homogeneous that the idea sometimes employed today of using a M\"obius strip topology to double the life of a belt made sense in general. Instead a M\"obius strip drive configuration was useful in the specific case of pulleys at right angles.

So might the idea of running a belt with a half-twist have been present long before the 1871 article in Scientific American? Leonardo da Vinci drew many belt-driven mechanisms in the 15th century, but we have not found a twisted belt among them. Where we have uncovered a twisted drive mechanism is in the golden age of Islamic science.

\section*{The M\"obius band chain pump of al-Jazari}

\begin{figure}[tbp]
\includegraphics*[height=0.6\textwidth]{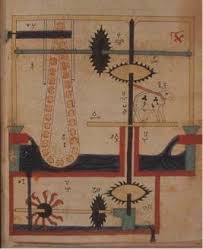}
\includegraphics*[height=0.6\textwidth]{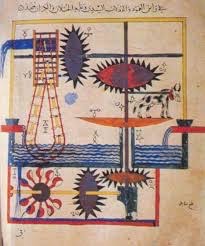} \\
\centering\includegraphics*[height=0.6\textwidth]{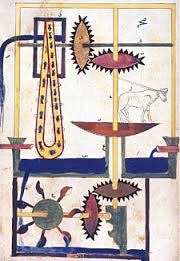}
\caption{\label{aljazari}
Chain pump showing M\"obius strip topology. Versions from three different copies of al-Jazari \cite{aljazari1206}.
}
\end{figure}

A fascinating illustration of an engineering mechanism displaying a M\"obius strip topology is found in the \emph{Book of Knowledge of Ingenious Mechanical Devices} written by al-Jazari in 1206 \cite{aljazari1206}. Al-Jazari, from Jazirat ibn Umar, present-day Cizre, Turkey, who introduced a number of novel mechanisms in this book, shows a M\"obius strip in the diagram of a chain pump (Fig.~\ref{aljazari}). Chain pumps have been utilized by many civilizations: the Babylonians, Egyptians, Chinese, and perhaps the Romans all used them before the Islamic Golden Age. However, this is the only case we have found in which the chain, or rope, linking the buckets for the pump is arranged as a M\"obius strip. In three versions of the diagram from three copies of the work that we show in Fig.~\ref{aljazari}  the pump drive has a M\"obius topology. 

\begin{figure}[tbp]
\centering\includegraphics*[height=0.6\textwidth]{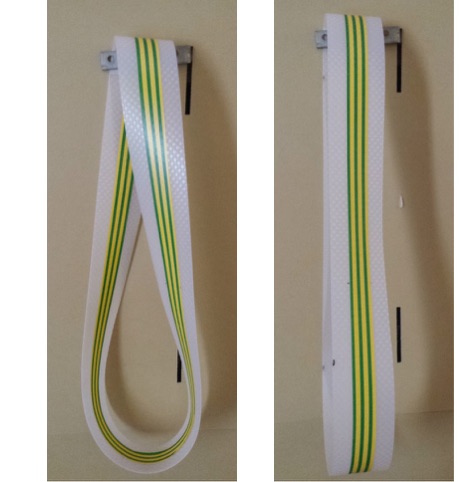}
\caption{\label{moebius1}
Juxtaposition of how a belt with normal two-sided topology hangs from a pulley to how a M\"obius strip hangs; compare with Figs~\ref{inman} and \ref{aljazari}.
}
\end{figure}

Recently Al-Hassani and Kiat  reconstructed al-Jazari's pump using computer assisted design \cite{alhassani2008}, but they did not take into account the M\"obius strip topology of the drive. It is clear that al-Jazari intended to represent a M\"obius strip. Following the edges of the belt, which in this case  is made with rope holding the water containers, shows that one has to travel twice around the circuit to close the path. Moreover, the work shows a refined use of perspective, showing a circle or cylinder as an ellipse. If 
a normal belt is hung from a pulley and observed perpendicular to the pulley one sees a rectangle as the front covers exactly the rear; however if one hangs a M\"obius strip the twist (note that this has chirality) allows us to observe both the front and back; Fig.~\ref{moebius1}. This is precisely what we see in Fig.~\ref{aljazari}. Of course this returns us to the quarter twist belt; a M\"obius strip is the natural solution precisely because it naturally hangs in this fashion.

The conclusion is that  al-Jazari accurately represented a M\"obius strip. From a mechanical point of view one may wonder how the containers could be arranged to pass inside the upper pulley, however the upper pulley is not solid, but a gear with horizontal bars, so the containers can pass between two bars and form a kind of toothed belt (another innovation). 
 In this way the containers are neither inside nor outside the belt, but pass through it; that is to say they form a `normal' belt while the M\"obius strip is determined by the rope, which is one long rope in this M\"obius strip configuration rather than the two ropes of half the length that would be needed  for a `normal' belt. Thus the containers are the same ones on one and the other side of the belt but pass once in one position and the next time rotated 180 degrees. An engineering advantage might be that containers could last longer being used symmetrically in this fashion and not always stressed on one side, and an asymmetrical breakage would still allow a damaged container to maintain a certain efficiency every second turn.

\section*{M\"obius strips before M\"obius?}

We consider that the Sentinum mosaic maker was quite aware of what he was doing: the  mosaic without a doubt depicts a M\"obius strip. As to whether it represents a M\"obius strip, that is, whether the artist was aware of the particular topological properties of the object depicted, that it has just one surface, probably not. The artist seems to have decided to depict a twisted strip in order to show all the zodiacal symbols, which could not be viewed from a single point of view if they were on one side of an untwisted belt; if one looks at the top of the hoop the artist emphasizes the edge and passes it from the bottom to the top when it would have been easier to join the upper and lower edges with continuity (we see in Fig.~\ref{aion3}(c)  a different, non-M\"obius  solution another artist employed). 
There are plenty of twisted bands with varying numbers of twists in Roman mosaics, and some of them are M\"obius strips, but we agree with Larison \cite{Larison1973}: we cannot say that the artists noted the topological implications. In the case of Al-Jazari it is clear that the belt has a single edge, i.e., a single rope, and the containers have to pass back and forth as they rotate. It is not clear that there was the notion that this forms an object with a single side; the interest is that the scheme should work in mechanical terms.
  As for the modern tradition of M\"obius strip drive belts, there we can see that by 1871 there was already an awareness of the topological implication of giving a belt a half-twist.  This American publication is intriguingly close in time, yet quite removed in space, from the `official' discoveries of M\"obius and Listing in Germany in 1858.  (Note that neither mathematician published in 1858; Listing published in 1861, while M\"obius never published. His work was found in an 1865 memoir and published posthumously after 1868  \cite{Fauvel1993,Pickover2006}.)  Had this knowledge crossed the Atlantic and diffused out into society to reach practical engineers in those 13 years, or has the modern engineering tradition  an independent, earlier origin? The latter seems likely, but that remains a conjecture until someone finds a pre-1858 reference to a half-twisted drive belt.
  The question of whether an understanding of the topological aspects of a M\"obius strip is an evolution of the knowledge and mathematical thinking developed thanks to the contributions of M\"obius and Listing remains open. The M\"obius strip predates the birth of modern topology; we may call it an \emph{ante literam} topology.

\section*{Acknowledgements}

JHEC acknowledges the financial support of the Spanish MINCINN project FIS2013-48444-C2-2-P.

\bibliographystyle{unsrt}
\bibliography{moebius}

\end{document}